\documentclass[twoside,leqno]{article}

\usepackage[letterpaper]{geometry}

\usepackage{siamproceedings}

\usepackage[T1]{fontenc}
\usepackage{amsfonts}
\usepackage{graphicx}
\usepackage{epstopdf}
\usepackage{enumitem}
\usepackage{algorithmic}
\ifpdf
  \DeclareGraphicsExtensions{.eps,.pdf,.png,.jpg}
\else
  \DeclareGraphicsExtensions{.eps}
\fi

\newsiamthm{rem}{Remark}
\newsiamremark{hypothesis}{Hypothesis}
\crefname{hypothesis}{Hypothesis}{Hypotheses}
\newsiamthm{claim}{Claim}
\newsiamthm{thm}{Theorem}
\newsiamthm{defi}{Definition}
\newsiamthm{que}{Question}
\newsiamthm{con}{Conjecture}
\newsiamthm{cor}{Corollary}
\newsiamthm{lem}{Lemma}
\newsiamthm{pro}{Proposition}
\newsiamthm{exe}{Example}

\usepackage{amsopn}

\newcommand{\La}{\Lambda}
\newcommand{\la}{\lambda}

\newcommand{\Corr}{{\rm Corr}}

\newcommand{\Spec}{{\rm Spec \,}}

\newcommand{\bif}{{\rm bif}}
\newcommand{\AS}{{\rm AS}}
\newcommand{\FS}{{\rm FS}}
\newcommand{\CE}{{\rm CE}}
\newcommand{\PCE}{{\rm PCE}}

\newcommand{\PR}{{\rm PR}}

\newcommand{\Gal}{{\rm Gal}}

\newcommand{\Rat}{{\rm Rat}}
\newcommand{\FL}{{\rm FL}}

\newcommand{\PGL}{{\rm PGL}\,}

\newcommand{\Aut}{{\rm Aut}}

\newcommand{\boxtensor}{{\Box\kern-9.03pt\raise1.42pt\hbox{$\times$}}}

\newcommand{\supp}{{\rm supp}\,}

\newcommand{\propsubset}
{\mbox{$\textstyle{
			\subseteq_{\kern-5pt\raise-1pt\hbox{\mbox{\tiny{$/$}}}}}$}}

\newcommand{\sH}{{\mathcal H}}

\newcommand{\sJ}{{\mathcal J}}

\newcommand{\sM}{{\mathcal M}}

\newcommand{\sO}{{\mathcal O}}
\newcommand{\sP}{{\mathcal P}}

\newcommand{\A}{{\mathbb A}}

\newcommand{\C}{{\mathbb C}}
\newcommand{\D}{{\mathbb D}}

\renewcommand{\P}{{\mathbb P}}
\newcommand{\Q}{{\mathbb Q}}
\newcommand{\R}{{\mathbb R}}

\newcommand{\Z}{{\mathbb Z}}

\newcommand{\bk}{{\mathbf{k}}}

\newcommand{\ep}{\varepsilon}

\newcommand{\Fix}{\mathrm{Fix}}

\begin{document}

\newcommand\relatedversion{}

\title{\Large Rigidity in Complex Dynamics: Multiplier Spectrum and Dynamical Andr\'e-Oort Conjecture\relatedversion}
    \author{Junyi Xie\thanks{Beijing International Center for Mathematical Research, Peking University (\email{xiejunyi@bicmr.pku.edu.cn}, \url{http://scholar.pku.edu.cn/xiejunyi}).}
    }

\date{}

\maketitle








\begin{abstract} 
	In this note, we present recent progress on rigidity problems in one-dimensional complex dynamics, including the proof of Dynamical Andr\'e-Oort conjecture for curves and generic injectivity of multiplier spectrum. The proofs combine ideas from algebraic geometry, Arakelov geometry and complex dynamics.
\end{abstract}

\section{Introduction.}
For $d\geq 2,$ let $\Rat_d(\C)$ be the space of degree $d$ rational maps on $\P^1(\C)$.  
 The group $\PGL_2(\C)= \Aut(\P^1(\C))$ acts on $\Rat_d(\C)$ by conjugacy. The geometric quotient 
$$\sM_d(\C):=\Rat_d(\C)/\PGL_2(\C)$$ is the (coarse) \emph{moduli space} of endomorphisms  of degree $d$ (c.f. \cite{Silverman2012}).
Let  $$\Psi: \Rat_d(\C)\to \sM_d(\C)$$ be the quotient morphism.  
By \cite[Theorem 4.36(c)]{Silverman2007},
the moduli space $\sM_d=\Spec (\sO(\Rat_d(\C)))^{\PGL_2(\C)}$ is an irreducible affine variety of dimension $2d-2$.
We may also view it as a complex orbifold (c.f. \cite{Milnor1993}, \cite{milnor2011dynamics}). Moreover, the moduli space is indeed defined over $\Z$ (c.f. \cite{Silverman2012}), in particular, for every algebraically closed field $\bk$, $\sM_d(\bk)$ is well-defined and parameterized the rational maps of degree $d$ on $\P^1$ defined over $\bk.$

\medskip

A rational map  is called \emph{Latt\`es} if it is semi-conjugate to an endomorphism on an elliptic curve. A Latt\`es map $f$ is called \emph{flexible Latt\`es} if one can continuously vary the complex structure of the elliptic curve to get a family of Latt\`es maps passing through $f$. The structure of flexible Latt\`es maps is well-understood \cite[Lemma 5.5]{milnor2006lattes}.
Let $\FL_d(\C)\subseteq \Rat_d(\C)$ be the locus of flexible Latt\`es maps, which is Zariski closed in $\Rat_d(\C)$.
Set 
$$\sM_d^{\star}(\C):=\sM_d(\C)\setminus \Psi(\FL_d(\C)).$$
We note that $ \FL_d(\C)=\emptyset$ when $d$ is not a square number,  and $ \Psi(\FL_d(\C))$ is an  algebraic curve when $d$ is a square number. 
We say that $f$ is of \emph{monomial type} if it is semi-conjugate to the map $z\mapsto z^n$ on $\P^1$ for some $|n|\geq 2.$ 
We call $f\in \Rat_d(\C)$ \emph{exceptional} if it is Latt\`es or of monomial type. One may check that $f$ is exceptional if and only if some iterate  $f^n$ is exceptional.
Exceptional endomorphisms are considered as the exceptional examples in complex dynamics.

\medskip
A fundamental problem in arithmetic/complex dynamics is to understand $\sM_d(\C)$ and its dynamically meaningful structures.

\subsection{Plan.}
In Section \ref{pcfdao}, I introduce notion of  PCF maps, which form an important subsets in the moduli space; the Dynamical Andr\'e-Oort (DAO) conjecture and its recent progress, in particular the solution of the DAO for curves by Ji and me.
In Section \ref{proofdao}, I sketch the proof of DAO for curves. In Section \ref{multiplierandlength}, I introduce the the notions of multiplier and length spectrum, which are fundamental invariants for rational maps. I discuss McMullen's multiplier spectrum rigidity theorem and the recent length spectrum rigidity theorem obtain by Ji and me. In Section \ref{genericinj}, I disuss the recent Theorem \ref{main} of Ji and me, which proves the generic injectivity of the multiplier spectrum morphism.


\section{PCF maps and DAO conjecture.}\label{pcfdao}
For any $f\in \Rat_d(\C)$ with $d\geq 2$, we denote by 
$$C(f):=\{x\in \P^1(\C)|\,\, df(x)=0\}$$ the set of critical points of $f$ and $$PC(f):=\cup_{n\geq 1}f^n(C(f))$$ the postcritical set. 
We say $f$ is {\em postcritically finite} (PCF) if its critical orbits are finite i.e. $\#PC(f)<\infty$. PCF maps play a fundamental role in complex dynamics
as the dynamical behavior of critical points usually reflect the general dynamical behavior of a rational map.
PCF maps are also interesting in arithmetic dynamics. In fact, as a consequence of Thurston's rigidity theorem \cite{Douady1993},  PCF maps are defined over $\overline{\Q}$ in the moduli space $\sM_d$, except for the well-understood one-parameter family of flexible Latt\`es maps.
Moreover, PCF maps are Zariski dense \cite[Theorem A]{de2018dynamical}. It was suggested by Silverman \cite{Silverman2012} that PCF points in the moduli space of rational maps $\sM_d$  play a role analogous to that played by CM points in Shimura varieties.

\medskip

One may consider PCF maps as "special points" in the moduli space. It is natural to ask what is the distribution of these special points.  
This leads to the Dynamical Andr\'e-Oort (DAO) conjecture proposed by Baker and DeMarco in \cite[Conjecture 1.10]{baker2013special} (see also \cite[Conjecture 1.1]{demarco2018critical} for the curve case i.e. $\dim \La=1$).

\begin{con}[DAO conjecture]\label{condao}
	Let $(f_t)_{t\in \La}$ be  a non-isotrivial algebraic family of rational maps with degree $d\geq 2$,  parametrized by an irreducible algebraic variety $\La$ over $\C$ of dimension $N$. Then the following are equivalent:
	\begin{enumerate}
		\item[(i)] $f_t$ is PCF for a Zariski dense subset of $t\in \La(\C)$;
		\item[(ii)] the family has at most $N$ independent critical orbit.
	\end{enumerate}
\end{con}
We explain the meaning of  ``independent critical orbit'' as follows:
Set $k:=\C(\La)$.
The geometric generic fiber of the family $(f_t)_{t\in \La}$ is a rational function $f_{\overline{k}}:\P^1_{\overline{k}}\to \P^1_{\overline{k}}$.  
Following DeMarco \cite{demarco2018critical},  a pair $a,b\in \P^1(\overline{k})$  is called {\em dynamically related } if there is an algebraic curve $V\subseteq \P^1_{\overline{k}}\times \P^1_{\overline{k}}$ such that  $(a, b)\in V$  and $V$ is preperiodic by the product map $f_{\overline{k}}\times f_{\overline{k}}:\P^1_{\overline{k}}\times \P^1_{\overline{k}}\to \P^1_{\overline{k}}\times \P^1_{\overline{k}}$.  
In other words,  a pair of points $a,b\in \P^1(\overline{k})$  are dynamically related if and only if the orbit of $(a,b)\in \P^1_{\overline{k}}\times \P^1_{\overline{k}}$ is not Zariski dense in $\P^1_{\overline{k}}\times \P^1_{\overline{k}}.$
We say that a family $(f_t)$ has at most $N$ independent critical orbit if for any $N+1$ critical points, there is a pair among them which is dynamically related. 
\medskip

One of the motivations of DAO conjecture comes from the analogy between PCF points and CM points and the Andr\'e-Oort Conjecture in arithmetic geometry, which was recently fully solved by  Pila, Shankar and Tsimerman \cite{pila2021canonical}.  DAO Conjecture and the Andr\'e-Oort Conjecture both fit the principle of unlikely intersections, however,  there is no overlap between them.

\medskip

In my work with Ji, we proved DAO conjecture when $N=1$ (c.f. \cite[Theorem 1.8]{ji2023dao}).
\begin{thm}[DAO for curves]\label{thmdao}
	Let $(f_t)_{t\in \La}$ be  a non-isotrivial algebraic family of rational maps with degree $d\geq 2$,  parametrized by an algebraic curve $\La$ over $\C$. Then the following are equivalent:
	\begin{enumerate}
		\item[(i)] $f_t$ is PCF for infinite $t\in \La(\C)$;
		\item[(ii)] the family has at most one independent critical orbit.
	\end{enumerate}
\end{thm}
We will sketch its proof in Section \ref{proofdao}.
Almost all previous results for DAO conjecture concern the curve case.
\subsection{Previous results of DAO For curves.}
In the fundamental work of Baker-DeMarco \cite{baker2013special}, they proved Conjecture \ref{condao} for families of polynomials parameterized by the affine line $\A^1(\C)$. Since Baker-DeMarco \cite{baker2013special}, plenty of works have been devoted to proving special cases of DAO Conjecture and its variations. 

Most progress is made for families of polynomials maps. See Ghioca-Hsia-Tucker \cite{ghioca2013preperiodic}, Ghioca-Krieger-Nguyen \cite{ghioca2016case}, Ghioca-Krieger-Nguyen-Ye \cite{ghioca2017dynamical}, Favre-Gauthier \cite{favre2018classification} \cite{favre2022arithmetic}, and Ghioca-Ye \cite{ghioca2018dynamical}.  Among these results, a remarkable work of Favre and Gauthier \cite{favre2022arithmetic} proved Theorem \ref{thmdao} in the 
polynomial case. 
\par For general rational maps, DeMarco-Wang-Ye \cite{de2015bifurcation} proved Theorem \ref{thmdao} for some dynamical meaningful algebraic curves in the moduli space of quadratic rational maps. Ghioca-Hsia-Tucker  \cite{ghioca2015preperiodic} proved a weak version of Theorem \ref{thmdao} for families of rational maps given by $f_t(z)=g(z)+t$, $t\in\C$, where $g\in \overline{\Q}(z)$ is of $\deg g\geq 3$ with a super-attracting fixed point at $\infty$.

\section{Proof of DAO for curves.}\label{proofdao}
In this section, we present a sketch of the proof of Theorem \ref{thmdao} (c.f. \cite{ji2023dao}). 
The direction that (ii) implies (i) was proved by DeMarco \cite[Section 6.4]{demarco2016bifurcations} using Montel's theorem. 
We only need to prove the direction that (i) implies (ii).
\subsection{Preparation.}
Let $d\geq 2$ and $\La$ be a smooth algebraic curve over $\C.$ Let 
\begin{align}\label{family}
	f:\La\times\P^1(\C)&\to \La\times\P^1(\C), \\
	(t,z)&\mapsto (t,f_t(z)),  \notag
\end{align}
be an endomorphism of degree $d$ over $\La.$ We call $f$ an \emph{algebraic family} of rational maps of degree $d$ over $\La.$
Moreover, we say that $f$ is an algebraic family over a subfield $K$ of $\C$ if both $\Lambda$ and $f: \La\times\P^1\to \La\times\P^1$ are defined over $K.$
In other words, give an algebraic family $f$ on a smooth algebraic curve $\La$ over $\C$ is  equivalent to give an algebraic morphism $\phi_f: \La\to \Rat_d, t\mapsto f_t.$
Moreover $f$ is defined over $K$ if $\La$ and $\phi_f$ are defined over $K.$ We say that $f$ is \emph{non-isotrivial} if $\Psi\circ\phi_f: \La\to \sM_d$ is not a constant map.

\medskip 
\par Let $\pi_1:\La\times\P^1(\C)\to \La$ and $\pi_2:\La\times\P^1(\C)\to \P^1(\C)$ be the canonical projections. Let $\omega_{\P^1}$ be the Fubini-Study form on $\P^1(\C)$, and let  $\omega_\La$ be a fixed K\"ahler form  on $\La(\C)$ with $\int_{\La} \omega_{\La}=1$. Let  $\omega_1:=\pi_1^\ast (\omega_{\La})$ and $\omega_2:=\pi_2^\ast (\omega_{\P^1})$.  The {\em relative Green current} of $f$ is defined by 
\begin{equation*}
	T_f:=\lim_{n\to+\infty} d^{-n} (f^n)^\ast (\omega_2)=\omega_2+dd^c g,
\end{equation*}
where $g$ is a H\"older continuous quasi-p.s.h. function \cite[Lemma 1.19]{dinh2010dynamics}. For every $t_0\in \La(\C)$, we have $$T_f\wedge [t=t_0]=\mu_{f_{t_0}},$$ where $\mu_{f_{t_0}}$ is the maximal entropy measure of $f_{t_0}$.
\par A {\em marked point} $a$ is a morphism $a:\La\to \P^1$. The {\em bifurcation measure} of the pair $(f,a)$ is defined by
\begin{equation*}
	\mu_{f,a}:=(\pi_1)_\ast (T_f\wedge [\Gamma_a])=a^*T_f,
\end{equation*}
where $\Gamma_a$ is the graph of $a$. 
By Gauthier-Vigny \cite[Proposition 13 (1)]{Gauthier2019}
$\mu_{f,a}$ has finite mass.

\medskip

\begin{defi}A marked point $a$ is called {\em active} if $\mu_{f,a}$ does not vanish. Otherwise, it is called {\em passive}. 
\end{defi}
A {\em marked critical point} $c$ is a marked point such that $c(t)$ is a critical point of $f_t$ for each $t\in \La(\C)$.  Let $f$ be an algebraic family of rational maps as in (\ref{family}) with marked critical points $(c_i)_{1\leq i\leq 2d-2}.$ The {\em bifurcation measure} of $f$ is defined by $$\mu_{\bif}:=\sum_{i=1}^{2d-2}\mu_{f,c_i},$$
which has finite mass. 
\medskip
The following theorem was proved by Dujardin-Favre \cite[Theorem 2.5]{Dujardin2008} for marked critical points and DeMarco \cite{demarco2016bifurcations} in general case.
\begin{thm}\label{stable}
	Let  $f$ be a non-isotrivial algebraic family, then $a$ is passive if and only if $a$ is  preperiodic.
\end{thm}

For the simplicity, we assume that the family $f: \La\times \P^1\to \La\times \P^1$ over $\La$ is {\bf non-isotrivial and defined over $\overline{\Q}$}.
Indeed we may reduce to such case using Thurston's rigidity theorem. After a suitable base change $\La'\to \La$, we may assume that {\bf we have exactly $2d-2$ marked critical points $c_1,\dots,c_{2d-2}$}. By contradiction, we may assume that {\bf $c_1,c_2$ are dynamically independent}.

\subsection{Overview:}
There are four steps in the proof.  As in the previous works, the first step is to show the equidistribution of PCF parameters. 
We then get an additional condition that for every $c_i$,  $\mu_{f,c_i}$ is proportional to the bifurcation measure $\mu_{\bif}$ and we only need to get a contradiction under this additional assumption. 
We do this in the next three steps using a new strategy. {\bf The basic idea is to work on general points with respect to $\mu_{\bif}$.}

In Step 2, we show  a selection  of  conditions are satisfied for $\mu_{\bif}$-a.e. point $t\in \La(\C)$, such that for a good parameter satisfies all these conditions, we can construct the similarities in Step 3, moreover  we can finally get a contradiction in Step 4.

In Step 3, we construct similarities between the phase space and the parameter space. Such similarities are known before in some cases for prerepelling parameters (which is a countable set). 
The novelty of this result is to get the similarities at $\mu_{\bif}$-a.e. point $t\in \La(\C)$. This leads to many difficulties coming  from the non-uniformly hyperbolic phenomenon.   

In Step 4, we  construct local symmetries of maximal entropy measures from the similarities constructed in Step 3,  and a contradiction comes from  these symmetries and an arithmetic condition that we selected in Step 2.

\subsection{Step 1: Equidistribution.}
The application of various equidistribution theorems is one of the most successful idea in arithmetic dynamics, which backs to the works of Ullmo \cite{Ullmo1998} and  Zhang \cite{Zhang1998}, in where they solved the Bogomolov Conjecture. It was first introduced to study Conjecture \ref{condao}
in Baker-DeMarco's fundamental work \cite{baker2013special}. 

In our setting, as $\La$ is not projective, we need the following result, which is a direct consequence of Yuan-Zhang's equidistribition theorem \cite[Theorem 6.2.3]{yuan2021} on quasi-projective varieties.
\begin{thm}\label{equi}
	Let $f$ be an algebraic family of rational maps as in (\ref{family}) and $a$ be  an active marked point, all defined over a number field $K$. Let $t_n\in \La(\overline{\Q})$ be an infinite sequence of distinct points such that $f_{t_n}$ is PCF. Then we have 
	$$\frac{1}{|\Gal(t_n)|}\sum_{t\in \Gal(t_n)}\delta_t\to\frac{\mu_{f,a}}{\mu_{f,a}(\La)},$$
	when $n\to+\infty$, where $\Gal(t_n)$ is the Galois orbit of $t_n$. 
\end{thm}

Yuan-Zhang's equidistribition theorem \cite[Theorem 6.2.3]{yuan2021} is based on their  recent theory of adelic line bundles on quasi-projective varieties.
Before Yuan-Zhang' theory, this step was usually non-trivial in the previous works as in \cite{baker2013special} and  \cite{favre2022arithmetic}.

The following result is a simple consequence of 
Theorem \ref{equi}.
\begin{cor}\cite[Corollary 2.4]{ji2023dao}\label{corequi}
	Let $f$ be an algebraic  family of rational maps as in (\ref{family}) with marked critical points $(c_i)_{1\leq i\leq 2d-2}$.  Assume that $f$ is defined over $\overline{\Q}$ and 
	there is an infinite sequence $t_n, n\geq 0$ of distinct points in $\La(\overline{\Q})$
	such that $f_{t_n}$ is PCF, then for $c_i, c_j$ being active marked critical points, we have
	$$\frac{\mu_{f,c_i}}{\mu_{f,c_i}(\La)}=\frac{\mu_{f,c_j}}{\mu_{f,c_j}(\La)}.$$
\end{cor}

Back to our setting, as condition (i) in Theorem \ref{thmdao} holds, Corollary \ref{corequi} shows that {\bf for every active $c_i$,  $\mu_{f,c_i}$ is proportional to the bifurcation measure $\mu_{\bif}$}.

\subsection{Step 2: Parameter exclusion.}

In this step, we show  a selection  of  conditions are satisfied for $\mu_{\bif}$-a.e. point $t\in \La(\C)$, such that for a good parameter satisfies these conditions, we can construct the similarities in Step 3, moreover,  we can finally get a contradiction in Step 4.

\medskip

Since aside from the flexible Latt\`es locus, the exceptional maps are isolated in the moduli space, it is easy to show that 
\begin{enumerate}
	\item[(1)] For $\mu_{\bif}$-a.e. point  $t\in \La(\C)$, $f_t$ is non-exceptional. 
\end{enumerate}

We say that $t\in \La(\C)$ is \emph{transcendental} if $t\not\in \La(\overline{\Q}).$ As $\mu_{\bif}$ does not have atom and $\La(\overline{\Q})$ is countable, we get
\begin{enumerate}
	\item[(2)]Let  $(f_t)_{t\in \La}$  be an algebraic family of rational maps over an algebraic curve $\La$,  then $\mu_{\bif}$-a.e. point  $t\in \La(\C)$ is transcendental. 
\end{enumerate}

\subsubsection{Invariant correspondences.}
Let
$\Corr(\P^1)^{f_t}_*$ be the set of $f_t\times f_t$-invariant Zariski closed subsets $\Gamma_t\subseteq \P^1\times \P^1$ of pure dimension $1$ such that both $\pi_1|_{\Gamma_t}$ and $\pi_2|_{\Gamma_t}$ are finite. Elements in $\Corr(\P^1)^{f_t}_*$ can be viewed as "dynamical self-relation" of $f_t$.

Let $\Corr^{\flat}(\P^1_{\La})^{f}_*$ be the set of $f\times_{\La}f$-invariant Zariski closed subsets $\Gamma\subseteq \La\times(\P^1\times \P^1)$ which is flat over $\La$ and whose generic fiber is in $\Corr(\P^1_{\eta})^{f_\eta}_*$, where $\eta$ is the generic point of $\La$. Elements in $\Corr^{\flat}(\P^1_{\La})^{f}_*$ can be viewed as "families of dynamical self-relation" of $f$.

In general, a correspondence $\Gamma_t\in \Corr(\P^1)^{f_t}_*$ may not be contained in any correspondence in $\Corr^{\flat}(\P^1_{\La})^{f}_*$.
On the other hand, it is the case if $t$ is transcendental by \cite[Proposition 3.11]{ji2023dao}.
Moreover, $\Corr^{\flat}(\P^1_{\La})^{f}_*$ is countable by \cite[Proposition 3.8]{ji2023dao}.

\medskip

Next we introduce the frequently separated condition.

\subsubsection{Frequently separated condition.}
Let $(M,d)$ be a metric space and $g: M\to M$ be a self-map.

\begin{defi}
	Let $A$ be a subset of $\Z_{\geq 0}$. The {\em asymptotic lower/upper density} of $A$ is defined by 
	\begin{equation*}
		\underline{d}(A):=\liminf_{n\to \infty} |A\cap[0,n-1]|/n,
	\end{equation*}
	and 
	\begin{equation*}
		\overline{d}(A):=\limsup_{n\to \infty} |A\cap[0,n-1]|/n.
	\end{equation*}
	If $\underline{d}(A)=\overline{d}(A)$, we set $d(A):=\underline{d}(A)=\overline{d}(A)$ and call it the  {\em asymptotic density} of $A$.
\end{defi}

\medskip

We still denote by $d$ the distance in on $M\times M$ by $$d((x_1,y_1),(x_2,y_2))=\max\{d(x_1,x_2),d(y_1,y_2)\}.$$
Let $\Sigma$ be a non-empty subset of $M\times M$. 
We view $\Sigma$ as a correspondence on $M.$ The most typical example is the diagonal. 
For $x\in M$, we denote by $\Sigma(x):=\pi_2(\pi_1^{-1}(x))$, where $\pi_1,\pi_2$ are the first and the second projections.
When $\pi^{-1}(x)\neq \emptyset$, for every $y\in M$, we have $d((x,y), \Sigma)\geq d(y, \Sigma(x)).$

\begin{defi}A pair of points $x,y\in M$ is called 
	\begin{enumerate}
		\item[(i)]{\bf Frequently separated} $\FS(\Sigma)$ for $\Sigma$, if for every $\ep>0$, there is $\delta>0$ such that 
		$$\overline{d}(\{n\geq 0|\, d((g^n(x),g^n(y)),\Sigma)\geq \delta\})>1-\ep.$$
		\item[(ii)]{\bf Average separated}  $\AS(\Sigma)$ for $\Sigma$, if 
		$$\liminf_{n\to \infty}\frac{1}{n}\sum_{i=0}^{n-1} \max\{-\log d((g^n(x),g^n(y)),\Sigma), 0\}<+\infty.$$
	\end{enumerate}
\end{defi}
The above conditions  depend only on the equivalence class of distance functions on $M$ and $M\times M.$
It is clear that if $\Sigma\subseteq \Sigma'$, then $\AS(\Sigma')$ (resp. $\FS(\Sigma')$) implies $\AS(\Sigma)$ (resp. $\FS(\Sigma)$).
Moreover we have 
\begin{lem}\cite[Lemma 3.5]{ji2023dao}\label{lemasimplyfs} The $\AS(\Sigma)$ condition implies the $\FS(\Sigma)$ condition.
\end{lem}

When $\Sigma$ is the diagonal $\Delta$, the $\FS(\Delta)$ condition means that in most of the time $n\geq 0$, the orbits of $x$ and $y$ are $\delta$-separated for some small $\delta>0$. The proportion of such time could tend to  $1$ when $\delta$ tends to $0.$

\medskip

Back to our setting, we may view $\P^1(\C)$ as a metric space. For every $\Gamma\in \Corr^{\flat}(\P^1_{\La})^{f}_*$, we have the conditions  $\AS(\Gamma_t)$, $\FS(\Gamma_t)$ with respect to the pair $c_1,c_2$ for every $t\in \La(\C)$.
Roughly speaking, $\FS(\Gamma_t)$ condition means that in most of the time $n\geq 0$, $(f_t^n(c_1(t)), f^n_t(c_2(t)))$ is not too close to $\Gamma.$
We have the following highly non-trivial statement 
\begin{enumerate}
	\item[(3)]Under the assumption (i) in Theorem \ref{thmdao}, $\mu_{\bif}$-a.e. point  $t\in \La(\C)$ satisfies  the  $\AS(\Gamma_t)$ (hence $\FS(\Gamma_t)$) condition for every $\Gamma\in \Corr^{\flat}(\P^1_{\La})^{f}_*$ (c.f. \cite[Corollary 3.14]{ji2023dao}). 
\end{enumerate}
To show that $\AS(\Gamma_t)$ is satisfied for $\mu_{\bif}$-a.e. point, we consider integrations with respect to $\mu_{\bif}$ having arithmetic meanings, and the aim is to show that these integrations  are bounded. We bound these integrations
via Yuan-Zhang's arithmetic intersection theory on quasi-projecitve varieties \cite{yuan2021}.

\subsubsection{Non-uniformly hyperbolic conditions.}
Let $a$ be any marked point. In this section, the distance and the norm of the derivatives are computed with respect to the metrics induced by $\omega_\La$ and $\omega_{\P^1}$. 
For every $n\geq 0$, let $\xi_{a,n}:\La\to \P^1(\C)$ denote the map $\xi_{a,n}(t):=f_t^n(a(t))$. 
\begin{defi}
	A parameter $t_0\in \La(\C)$ is called 
	\begin{enumerate}
		\item[(i)] {\bf Marked Collet-Eckmann} $\CE^*(\la)$ for some $\la>1$, if there exists $C>0$ and $N>0$ such that $$|df_{t_0}^n(f_{t_0}^N(a(t_0)))|\geq C\la^n$$ for every $n\geq 0$. 
		\item[(ii)]  {\bf Parametric Collet-Eckmann} $\PCE(\la)$ for some $\la>1$, if there exists $C>0$ such that $$\left|\frac{d\xi_{a,n}}{dt}(t_0)\right|\geq C\la^n$$ for every $n\geq 0$. 
		\item[(iii)]  {\bf Polynomial Recurrence} $\PR(s)$ for some $s>0$, if there exists an integer $N>0$ such that
		$$ d(f_{t_0}^n(a(t_0)), \mathcal{C}_{t_0})\geq n^{-s}$$
		for every $n\geq N$.  Where $\mathcal{C}_{t_0}$ is the critical set of $f_{t_0}$ 
	\end{enumerate}
\end{defi}

A result of  De Th\'elin-Gauthier-Vigny \cite{de2021parametric} shows that 
\begin{enumerate}
	\item[(4)]There is $\la>1$ such that $\mu_{f,a}$-a.e. parameters satisfy the $\CE^*(\la)$ and the $\PCE(\la)$ for 
	every active marked critical point $c_i$.
\end{enumerate}

Using pluripotential theory, we proved the following
\begin{enumerate}
	\item[(5)]For every compact subset $K\subset \La(\C)$, there exists $s>0$ depending on $K$,  such that  $\mu_{f,a}$-a.e. $t\in K$  satisfy the $\PR(s)$ condition (c.f.\cite[Theorem 4.6]{ji2023dao}).
\end{enumerate}

\begin{defi}\label{ce}
	A rational map $g$ of degree at least $2$ is called {\em Collet-Eckmann} $\CE(\la)$ for some $\la>1$ if:
	\begin{enumerate}
		\item[(i)] There exists $C>0$ such that  for every critical point $c\in\sJ(g)$, there exists $N>0$ such that $|dg^n(g^N(c))|\geq C\la^n$ for every $n\geq 1$;
		\item[(ii)] $g$ has no parabolic cycle. 
	\end{enumerate}
\end{defi}

\par The proof of Condition (6)  is a combination of Condition (4),  Siegel's linearization theorem \cite[Theorem 11.4]{milnor2011dynamics} and the fact that 
the set of Liouville numbers has Hausdorff dimension $0$  \cite[Lemma C.7]{milnor2011dynamics},
\begin{enumerate}
	\item[(6)]  Let  $(f_t)_{t\in \La}$  be an algebraic family of rational maps over an algebraic curve $\La$, assume moreover that  for every active  marked critical point $c_i$,  $\mu_{f,c_i}$ is proportional to the bifurcation measure $\mu_{\bif}$, then  $\mu_{\bif}$-a.e. point $t\in \La(\C)$ satisfies the  Collet-Eckmann condition.
\end{enumerate}
\par Condition (1), (2) and (6) are relatively easy to show, and Condition (4) can be easily deduced by the work of 
De Th\'elin-Gauthier-Vigny \cite{de2021parametric}. The major part of Step 2 is the proof of Condition (3) and (5). The proof of Condition (5) requires pluripotential theory. The proof of Condition (3) requires both pluripotential theory and arithmetic intersection theory.

\subsection{Step 3: Similarity between the phase space and the parameter space.}
To get a contradiction, we show that the conditions (1),(2),(4),(5),(6) imply the opposite of (3).
Our idea is to get similarity between the bifurcation measure $\mu_{\bif}$ on the parameter space and the maximal entropy measure $\mu_{f_t}$ on the phase space.
This can be thought  of as a generalization of Tan's work \cite{Tan1990}, in where she got such similarities at Misiurewicz points in Mandelbrot set,  and as a generalization of Gauthier's \cite[Section 3.1]{gauthier2018dynamical} and Favre-Gauthier's works \cite[Section 4.1.4]{favre2022arithmetic}, in where they got such similarities at properly prerepelling parameters.
\medskip

The key is to construct renormalization maps.
\subsubsection{Renormalization maps.}
Recall that $\xi_{a,n}:\La\to \P^1$ is the map $\xi_{a,n}(t):=f_t^n(a(t))$. 
\begin{thm}\cite[Theorem 5.2]{ji2023dao}\label{thmsimlarity}
	For $t_0\in \La(\C)$ satisfies (4),(5),(6), let $a$ be any active marked point.
	there is a subset $A\subseteq \Z_{\geq 0}$ of large lower density and a sequence of positive real numbers $(\rho_n)_{n\in A}$ tending to zero, such that we can construct a 
	family of renormalization maps 
	\begin{align*}
		h_n:\D&\to \P^1(\C), n\in A
		\\t&\mapsto \xi_{a,n}\left( \rho_n t\right)
	\end{align*}
	where the unit disk $\D:=\{z\in \C|\,\, |z|<1\}$ is identified with a small disk centered at $t_0$. We show that this family is normal and no subsequences of $h_n, n\in A$ can tends to a constant map.
\end{thm}

For our application, we take $a$ to be any active marked critical point $c_i$.
\begin{rem}In fact, for this result, we may replace the Collet-Eckamnn (=CE) condition by the \emph{Topological Collet-Eckamnn condition (=TCE)}. By Przytycki-Rohde \cite{przytycki1998porosity} (see also \cite[Main Theorem]{przytycki2003equivalence}), the CE condition implies the TCE condition.
\end{rem}

Comparing with the previous results, we get similarities between phase space and parameter space not only for prerepelling parameters (which is a countable set) but also for parameters satisfying Topological Collet-Eckamnn condition and Polynomial Recurrence condition  (which is a set of full $\mu_{\bif}$ measure).  
I believe that this result has an independent interest in complex dynamics.

In the previous works for prerepelling parameters, the key point is that the orbits of the marked points have a uniform distance from the critical locus for all but finitely many terms.
This is not true in our case. For this reason, we introduce the following new strategy. 

\subsubsection{Proof strategy for Theorem \ref{thmsimlarity}.}
The proof is divided into two parts.

In the first part, we work only on the phase space. 
We select the ``good time set" $A\subseteq \Z_{\geq 0}$ of large lower density. For each good time $n\in A$, we construct $n$ maps from certain fixed simply connected domain to $\P^1(\C)$. 
Roughly speaking, the goodness of $n$ means that the above maps have a uniformly bounded number of critical points. Then we study the distortions of such maps which are  non-injective in general (c.f. \cite[Section 6]{ji2023dao}).  To describe  the distortion of  perhaps non-injective holomorphic maps, we introduce the concepts of upper and proper lower radius (c.f. \cite[Section 5]{ji2023dao}).
Comparing with the usual lower radius of the image, the advantage of the proper one is the stability under small perturbations.

In the second part, we use a binding argument to get the renormalization maps from the maps defined above and show that this family is normal and no subsequence of $h_n, n\in A$ can tends to a constant map  (c.f. \cite[Section 7]{ji2023dao}). In particular, we decide the rescaling factors $\rho_n, n\in A$ in this process.

%

\subsubsection{Similarity.} Identify $t_0$ in Theorem \ref{thmsimlarity} with $0$ in $\D.$ As in Theorem \ref{thmsimlarity}, $a$ is an active marked point, e.g. any active marked critical point $c_i.$

Let $\mu_{f_0}$ be the maximal entropy measure of $f_0$ ($=f_{t_0}$ under the identification).
For $\rho\in (0,1]$, let $[\rho]: \D\to \D$ be the map $t\mapsto \rho t$. 
The following result shows that under the assumption of Theorem \ref{thmsimlarity},
$\mu_{f_0}$ can be read from $\mu_{f,a}$ using the renormalization maps $h_n, n\in A$.


\begin{pro}\cite[Proposition 7.1]{ji2023dao}\label{measure}
	Let $f:\D\times \P^1\to \D\times \P^1$ be a holomorphic family of rational maps as in (\ref{family}) and $a$ be a marked point. 
	Let $n_j, j\geq 0$ be an infinite subsequence of $n\geq 0.$
	Let $\rho_{n_j}, j\geq 0$ be a sequence of positive real number tending to $0.$
	Define $h_{n_j}:=\xi_{a,n_j}\circ [\rho_{n_j}].$
	Assume that $h_{n_j}\to h$ locally uniformly. Then we have 
	\begin{equation}\label{equmealimju}d^{n_j}[\rho_{n_j}]^*\mu_{f,a}\to h^\ast (\mu_{f_0})
	\end{equation}
	where the convergence is the weak convergence of measures. 
	Moreover, if $a(0)\in \sJ(f_0)$ and $h$ is non-constant, then $0\in \supp\,h^\ast (\mu_{f_0}).$
\end{pro}

\begin{defi} We say that two sequences of positive real numbers $\rho_n, n\geq 0$ and $\rho_n', n\geq 0$ are \emph{equivalent} if $\log (\rho_n/\rho'_n), n\geq 0$ is bounded.
\end{defi}

The next result shows that the scales $\rho_n, n\in A$ are determined by $\mu_{f,a}$ itself up to equivalence. Hence $\mu_{f_0}$ can be read from $\mu_{\bif}$.
\begin{pro}\cite[Proposition 7.2]{ji2023dao}\label{proresunique}
	Let $\mu$ be a Borel  measure on $\D$.
	Let $\rho_n,\rho'_n, n\geq 0$ be two sequences of real numbers in $(0,1]$. Let $d_n, n\geq 0$ be a sequence of positive real numbers.
	Let $\mu_1,\mu_2$ be Borel  measures on $\D$ having positive mass.
	Assume that $\mu_1(\{0\})=\mu_2(\{0\})=0$.
	If $d_n[\rho_n]^*\mu\to \mu_1$ and $d_n[\rho_n']^*\mu\to \mu_2$, then 
	$\rho_{n_j}, j\geq 0$ and $\rho'_{n_j}, j\geq 0$ are equivalent.  
\end{pro}
\subsection{Step 4:  Conclusion via local symmetries of maximal entropy measures.}
\par We have constructed a family of renormalization maps $h_n:\D\to \P^1, n\in A$.  Since $A$ has a large lower density, after taking an intersection, we may assume that the set $A$ for $c_1$ and $c_2$ are the same.
We denote by $h_{c_1,n}, n\in A$ and $h_{c_2,n}, n\in A$ the renormalization maps for $c_1$ and $c_2$ respectively.

After suitable adjustments of $\sH_{c_1}:=\{h_{c_1,n}, n\in A\}$ and $\sH_{c_2}:=\{h_{c_2,n}, n\in A\}$, we show that they form an \emph{asymptotic symmetry} (c.f. \cite[Section 7.1]{ji2023dao}), which basically means that 
every limit of $h_{a,n}\times h_{b,n}$ produces a holomorphic local symmetry of $\mu_{f_t}$. We use such local symmetries and a compactness argument to show that condition (3) is not true (c.f. \cite[Proposition 7.8]{ji2023dao}). This gives a contradiction.
A key ingredient to prove \cite[Proposition 7.8]{ji2023dao} is \cite[Theorem 1.7]{jixielocal2022} of Ji and me, which shows that 
every holomorphic local symmetry of the maximal entropy measure comes from an (algebraic) invariant correspondence.
Since we may assume that $f_t$ is Collet-Eckmann, in this last step we may replace \cite[Theorem 1.7]{jixielocal2022} by 
combining \cite[Theorem A]{dujardin2022two} with \cite[Corollary 3.2]{dujardin2022two}.

\section{Multiplier spectrum and length spectrum.}\label{multiplierandlength}
As $\sM_d$ is affine, it is natural to study its ring of functions.
There is  a natural  dynamically meaningful family of morphisms from $\sM_d(\C)$ to  affine spaces, which we call the \emph{multiplier spectrum morphisms}. 
They give a natural system of 
functions on $\sM_d(\C)$, which can be viewed as analogies of the theta functions on the moduli space of elliptic curves (c.f. \cite[Remark 5.9]{Silverman2023}).

\medskip

For every $f\in \Rat_d(\C)$ and $n\geq 1$, denote by $\Fix(f^n)$ the multi-set of $f^n$-fixed points. Then $\Fix(f^n)$ has 
exactly $N_n:=d^n+1$ fixed points counted with multiplicity. 
For any $f^n$-fixed point $x$, the differential $df^n(x)\in \C$ is called the
the \emph{multiplier} of $f^n$ at $x$. 
Define $$S_n(f):=\{df^n(x)|\,\, x\in \Fix(f^n)\}\in \C^{N_n}/\Sigma_{N_n}$$
where $\Sigma_{N_n}$ is the symmetric group which acts
on $\C^{N_n}$ by permuting the coordinates. Via elementary symmetric polynomials, we have $\C^{N_n}/\Sigma_{N_n}\simeq \C^{N_n}.$
\begin{defi}
The \emph{multiplier spectrum} of $f$ is the sequence $S_n(f), n\geq 1.$
\end{defi}
According to the above definition, we need to collect the multipliers of all periodic points. This places the values of the multiplier spectrum in an infinite dimensional space, which could be inconvenient sometimes. However, by the Noetherianity, we can truncate it to certain period as follows, obtaining an algebraic morphism $\tau_d$ with values in a finite dimensional space while preserving most of the information of the multiplier spectrum. The precise definition is as follows:
Since $S_n$ takes the same value in a conjugacy class of rational maps, it descents to a morphism $$S_n:\sM_d(\C)\to \C^{N_n}/\Sigma_{N_n}.$$
For every $n\geq 1$, define $\tau_{d,n}$  be the morphism
\begin{align*}
	\tau_{d,n}:\sM_d(\C)&\to \C^{N_1}/\Sigma_{N_1}\times\cdots\times \C^{N_n}/\Sigma_{N_n},\\
	[f]&\mapsto (S_1(f),\dots,S_n(f)).
\end{align*}
For every $f\in \Rat_d(\C)$, set $\tau_{d,n}(f):=\tau_{d,n}([f])$.
For each $n\geq 1$, set $$R_n:=\{([f],[g])\in \sM_d(\C)^2|\,\,\tau_{d,n}([f])=\tau_{d,n}([g])\}.$$
Then  $R_n, n\geq 1$  form a decreasing sequence of Zariski closed
subsets  of $ \sM_d(\C)^2$.
By the Noetherianity, the sequence $R_n$ is stable for $n$ sufficiently large. Hence there exists a minimal positive  integer $m_d$ such that  $\tau_{d,m_d}(f)=\tau_{d,m_d}(g)$ implies that  $\tau_{d,n}(f)=\tau_{d,n}(g)$ for every $n\geq 1$, i.e. $f$ and $g$ have the same multiplier spectrum.  We define $$\tau_d:= \tau_{d,m_d}.$$
\begin{que}What is the exact value of $m_d$?
	\end{que}
It is also interesting to know its lower/upper bounds.


\medskip

Replace the multipliers by its norm in the definition of multiplier spectrum, and one get the definition
of the length spectrum. More precisely, for every  $f\in \Rat_d(\C)$ and $n\geq 1,$ denote by
$$L_n(f):=\{|df^n(x)|x\in\Fix(f^n)\}\in \R^{N_n}/\Sigma_{N_n}.$$

\begin{defi}The \emph{length spectrum} of $f$ is defined to be the sequence $L_n(f), n\geq 1$.
	\end{defi}
A priori, the length spectrum contains less information than the multiplier spectrum. 
But in \cite[Conjecture 1.10]{Ji2023b}, Ji and me we proposed the following description for rational maps with the same length spectrum.
\begin{con}\label{1.9}
	Let $f,g$ be non-Latt\`es rational maps of degree $d\geq 2$. If $f$ and $g$ has the same length spectrum, then $\tau_d(f)$ equals to either $\tau_d(g)$ or $\tau_d(\overline{g})$. Here $\overline{g}$ is the complex conjugation of $g$.
\end{con}

\medskip

In 1966, Kac  asked the famous question (c.f. \cite{Kac1966}):
\begin{center}{\em "Can one
		hear the shape of a drum?"}
\end{center}
Kac's question has been investigated in various settings in dynamics and geometry. In the contexts of geodesic flows on Riemannian manifolds with negative curvature,  
a long-standing conjecture stated by Burns-Katok \cite{Burns1985} (and
probably be considered even before) asserted
the rigidity of marked length spectrum (for closed geodesics).  The surface case was proved by Otal \cite{Otal1990} and by Croke \cite{Croke1990} independently. 
A Local version of the Burns-Katok conjecture in any dimension
was proved by Guillarmou-Lefeuvre \cite{Guillarmou2019}. 
In one-dimensional real dynamics, it was proved for expanding circle maps (see Shub-Sullivan \cite{shub1985expanding}), and for some unimodal maps (see Martens-de Melo \cite{martens1999multipliers} and Li-Shen \cite{li2006smooth}).
It was also studied in dynamical billiards.  We refer the readers to Huang-Kaloshin-Sorrentino \cite{Huang2018}, B\'alint-De Simoi-Kaloshin-Leguil \cite{Balint2020}, De Simoi-Kaloshin-Leguil \cite{DeSimoi2023},  and the references therein.

\medskip

The following question can be viewed as Kac's question in one-dimensional complex dynamics.
\begin{que}\label{questionmullendetf}For $f\in \Rat_d$, how does $S(f)$ (resp. $L(f)$) determine $f$?
\end{que}

\subsection{Rigidity of multiplier/length spectrum.}
The first remarkable answer Question \ref{questionmullendetf} is the following result proved by McMullen \cite{McMullen1987} in 1987.
\begin{thm}[Multiplier spectrum rigidity]\label{thmmcmullen}
Aside from the flexible Latt\`es family, the multiplier spectrum determines the conjugacy class of $f\in\Rat_d(\C)$, $d\geq 2$, up to finitely many choices.
\end{thm}
McMullen's proof relies on Thurston's rigidity theorem for PCF maps \cite{Douady1993}, in where Teichm\"{u}ller theory is essentially used. 
In my paper with Ji \cite{ji2023homoclinic}, we gave a new proof of McMullen's theorem without using  quasiconformal maps or Teichm\"{u}ller theory.
A key of the proof is to prove a criterion for $f\in \Rat_d(\C)$ being
exceptional via the information of a homoclinic orbit (c.f. \cite[Theorem 2.11]{ji2023homoclinic}),
whose proof relies only on some basic complex analysis.
Except \cite[Theorem 2.11]{ji2023homoclinic}, our proof of Theorem \ref{thmmcmullen} only requires some basic knowledges in Berkovich dynamics and hyperbolic dynamics.

The statement of Theorem \ref{thmmcmullen} is purely algebraic. However, there is no purely algebraic proof so far. It will be very interesting to have such a proof, as it may implies certain generalization of positive characteristic.

\medskip

Ji and I proved a parallel result for length spectrum \cite[Theorem 1.5]{ji2023homoclinic}, which 
strengthens McMullen's rigidity theorem.
\begin{thm}[Length spectrum rigidity]\label{thmlength}
	Aside from the flexible Latt\`es family, the length spectrum determines the conjugacy class of $f\in \Rat_d(\C)$, $d\geq 2$, up to finitely many choices.
\end{thm}

The multiplier spectrum (hence the length spectrum) is constant in any irreducible component of $\FL_d(\C)$, so it is necessary to exclude such family in Theorem \ref{thmmcmullen} and Theorem \ref{thmlength}.

\medskip

The length spectrum is harder to treat than the multiplier spectrum. A major difficulty is that the length spectrum map is not real algebraic. To be more precise, we view $\text{Rat}_d(\C)$ as a real algebraic variety by splitting the complex variable to two real variables via $z=x+iy$. A more theoretical way to do this is using the notion of Weil restriction (c.f. \cite[Section 4.6]{Poonen2017}). It is not hard to see that for every $\alpha\in \R^{N_n}/\Sigma_{N_n}$, the preimage $L_n^{-1}(\alpha)$ under the $n$-th Length spectrum map is the image of a closed real algebraic subset in some real algebraic variety under a finite morphism. However, unlike the complex case, such an image is merely a closed semi-algebraic set, but not real algebraic in general. 
\begin{exe}The images of $\A^1(\R)$ under the finite morphisms $z\mapsto z^2$ is $[0,+\infty)\subseteq \A^1(\R)$ which is not real algebraic.
	\end{exe}
In fact, in my paper with Ji \cite{ji2023homoclinic}, we gave an explicit example that $L_n^{-1}(\alpha)$ is not real algebraic:
\begin{exe}\cite[Theorem 8.10]{ji2023homoclinic}\label{thmexamplenonalg}
	Consider the first length spectrum map $L_1: {\rm Rat}_2(\C)\to \R^3/\Sigma_3.$
	For every $a\in (1,\sqrt{2})$, $\alpha:=\{a,a,a\}\in \R^3/\Sigma_3$, $L_1^{-1}(\alpha)$ is not real algebraic in ${\Rat}_2(\C).$
\end{exe}
This cause a problem as closed semi-algebraic sets do not satisfy the descending chain condition. 
\begin{exe}
For $n\in \Z_{\geq 0}$, set $Z_n:=[n,\infty)\subseteq \A^1(\R).$ We have $Z_{n+1}\subset Z_n$ but $\cap_{n\geq 0} Z_n=\emptyset$. 
\end{exe}
To solve this problem, our idea is to introduce a class of closed semialgebraic sets called  \emph{admissible} subsets (c.f. \cite[Section 8.2]{{ji2023homoclinic}}).
Roughly speaking, admissible subsets are the closed subsets which are images of real algebraic subsets under \'etale morphisms. We show that admissible subsets satisfy the descending chain condition. We showed that for every $\alpha=\{a_1,\dots,a_{N_n}\}\in (\R\setminus \{1\})^{N_n}/S_{N_n}$, $L_n^{-1}(\alpha)$ is admissible (c.f. \cite[Proposition 8.21]{{ji2023homoclinic}}). This is sufficient to prove Theorem \ref{thmlength}, because except finitely many of them, all periodic points are repelling. 

\medskip 
 
To prove Theorem \ref{thmlength}, except the notion of admissible subsets above, there are two key ingredients: One is Sullivan's rigidity theorem \cite{sullivan1986quasiconformal} for non-linear \emph{conformal expending repeller (=CER)}. Another is our \cite[Theorem 1.1]{ji2023homoclinic}, which shows that the existence of linear CER implies the exceptionalness.

\subsubsection{Uniform rigidity.}  
McMullen's multiplier spectrum rigidity Theorem (i.e. Theorem \ref{thmmcmullen}) equivalents to say that the morphism  $\tau_d$ is quasi-finite. So 
 there is $N\geq 1$ depending only on $d\geq 2$, such that for every $f\in \text{Rat}_d(\C)\setminus FL_d(\C)$, 
$$\#\Psi(\{g\in \text{Rat}_d(\C)\setminus FL_d(\C)|\,\, S_i(g)=S_i(f), i=1,\dots,N\})\leq N.$$
This can be viewed as a uniform version of Theorem \ref{thmmcmullen}.

In \cite[Question 1.6]{ji2023homoclinic}, Ji and I asked  whether the uniform version of length spectrum rigidity Theorem (i.e. Theorem \ref{thmlength}) holds.
\begin{que}\label{queuniformlen}
	Is there $N\geq 1$ depending only on $d\geq 2$, such that for every $f\in \text{Rat}_d(\C)\setminus FL_d(\C)$, 
	$$\#\Psi(\{g\in \text{Rat}_d(\C)\setminus FL_d(\C)|\,\, L_i(g)=L_i(f), i=1,\dots,N\})\leq N\,?$$
\end{que}

For every $n\geq 0,$ we set $$R_n:=\{(f,g)\in (\text{Rat}_d(\C)\setminus FL_d(\C))^2|\,\, L_i(f)=L_i(g), i=1,\dots, n\},$$
It is a decreasing sequence of closed subsets of $(\text{Rat}_d(\C)\setminus FL_d(\C))^2.$ 
If the sequence 
$R_n, n\geq 0$ is stable, for example if one can show that $R_n$ are admissible, 
then Question \ref{queuniformlen} has a positive answer. 
But at the moment, we only know that $R_n$ are semialgebraic.

\section{Generic injectivity of multiplier spectrum.}\label{genericinj}
		\begin{defi}
	For a point $x\in\sM_d(\C)$, we say that $\tau_d$ is injective at $x$ if $\tau_d^{-1}(\tau_d(x))=\left\{x\right\}$.
	For  a subset $X\subset \sM_d(\C)$, we say that $\tau_d$ is injective on $X$ if $\tau_d$ is injective at every $x\in X$.
\end{defi}
\par We quote the following question about the injectivity of $\tau_d$  from  McMullen \cite[Page 489]{McMullen1987}:
\medskip
\par {\em Noetherian properties imply there are an $N$ and an $M$ such that $E_1(R),\dots,E_N(R)$  determine $R$ up to at most $M$ choices, $\dots$ is $R$ determined uniquely?}\footnote{Following McMullen's notations, $E_i(R)$ is the multiplier spectrum of periods $i$ of the rational map $R$.}
\medskip


It turns out that $\tau_d$  is {\bf not} always  injective on $\sM_d^{\star}(\C)$.  
Silverman gave such examples by considering rigid Latt\`es map \cite[Theorem 6.62]{Silverman2007}. 
Another construction was introduced by Pakovich \cite{Pakovich2019}: Let $f$ be a rational map. Following Pakovich \cite{Pakovich2019}, for any decomposition $f=h_1\circ h_2$ into a composition of rational maps, we say that the rational map $\tilde{f}:=h_2\circ h_1$ is an \emph{elementary transformation} of $f$.  \begin{defi}We say that rational maps $f$ and $g$ are \emph{elementarily equivalent} if there exists a chain of elementary
transformations between $f$ and $g$. 
\end{defi}
Easy to check that if $f$ and $g$ are elementarily equivalent then $\tau_d(f)=\tau_d(g)$
(c.f. \cite[Lemma 2.1]{pakovich2019recomposing}).  

\medskip

Even though  $\tau_d$  is not always injective on $\sM_d^{\star}(\C)$ as we have seen, one might hope that $\tau_d$ is  \emph{generically injective}, i.e. $\tau_d$ is injective on a Zariski open subset. 
Poonen asked whether $\tau_d$ is always generically injective \cite[Question 2.43]{Silverman2012}.  
In my paper with Ji \cite{Ji2023b}, we proved the following result.
It  gives an affirmative answer to Poonen's question, hence an affirmative answer to McMullen's question for generic parameters. 
\begin{thm}\cite[Theorem 1.3]{Ji2023b}\label{main}
	For every $d\geq 2$,
	the morphism $$\tau_d:\sM_d(\C)\to \C^{N_1}/\Sigma_{N_1}\times\cdots\times \C^{N_{m_d}}/\Sigma_{N_{m_d}}$$ is generically injective. 
\end{thm}
As a consequence, the system of functions given by the multiplier spectrum generate the function field of $\sM_d.$  

\subsection{Sketch of the proof of Theorem \ref{main}.}
We first introduce two key ingredients in the proof.
\begin{defi}\label{defiintert}
	Let $d\geq 2$ and $f,g\in\Rat_d(\C)$. We say that $f$ and $g$ are \emph{intertwined} if there exists a (maybe reducible) algebraic curve $Z\subset \P^1\times \P^1$ whose projections to both factors are finite surjective,
	and $Z$ is invariant by the map $f\times g:\P^1\times \P^1\to \P^1\times \P^1$.  
\end{defi}

Our first ingredient is the following DAO type result.

\begin{thm}\cite[Theorem 1.5 and Remark 1.6]{Ji2023b}\label{DAOint}
	Let $d\geq 2$ and let $f_V,g_V$ be two degree $d$ non-isotrivial algebraic families defined over $\C$, parametrized by the same irreducible algebraic curve $V$, and $f_V,g_V$  are not families of flexible Latt\`es maps. Assume that there are infinitely many $t\in V(\C)$ such that $f_t$ and $g_t$ are both PCF. Then for every $t\in V(\C)$, $f_t$ and $g_t$ are intertwined. 
\end{thm}
Our proof for Theorem \ref{DAOint} is similar to the proof of Theorem \ref{thmdao} (DAO for curves).

\medskip

The second ingredient is the following theorem of Pakovich.
\begin{thm}\cite[Theorem 1.2]{Pakovich2025}\label{intertwined}
	Let $d\geq 2$. Then there exists a Zariski dense open subset $W_d$ of $\Rat_d(\C)$ such that for every $f,g\in W_d$, they are intertwined if and only if they are conjugate to each other. 
\end{thm}
The above result \cite[Theorem 1.2]{Pakovich2025} was proved after our paper \cite{Ji2023b}. The in the original proof of Theorem \ref{main} in \cite{Ji2023b}, we use a previous paper \cite{pakovich2021iterates} of Pakovich, where he established Theorem \ref{intertwined} for $d\geq 4.$
This is sufficient to prove Theorem \ref{main} for $d\geq 4.$ However, the cases $d=2,3$ was known from earlier works of Milnor \cite{Milnor1993} ($d=2$ case) and  Gotou \cite[Theorem 1.2]{Gotou2023} ($d=3$ case).

\medskip

By contradiction, assume  that $\tau_d$ is not generically injective. 
We consider the following special class for PCF points in $\Rat_d(\C).$
\begin{defi}
	A PCF map $f\in \Rat_d(\C)$ is called \emph{hyperbolic of disjoint type} if $f$ has  $2d-2$ distinct super-attracting cycles.  
\end{defi}
By \cite[Theorem 6.6]{gauthier2019hyperbolic}, such points are Zariski dense in $\Rat_d(\C)$.
We construct two algebraic families of rational maps $f_{t}$ and $g_{t}$, parametrized by the same algebraic curve $V$, such that: 
\begin{enumerate}
	\item[(1)] for every $t\in V(\C)$, the images $f_{t}$ and $g_{t}$ in $\sM_d$ are different; 
	\item[(2)]  $\tau_d(f_{t})=\tau_d(g_{t})$ for every $t\in V(\C)$;
	\item[(3)] for every $t\in V(\C)$, $f_t$ and $g_t$ are contained in $W_d$;
	\item[(4)] there are infinite $t\in V(\C)$ such that $f_t$ is hyperbolic PCF of disjoint type.
\end{enumerate}
An elementary and combinatorial argument show that whether a rational map $f$ is hyperbolic PCF can be read from its multiplier spectrum $S(f)$ (c.f. \cite[Lemma 3.5]{Ji2023b}).
By (2) and (4) above, there are infinitely many $t\in V(\C)$ such that both $f_t$ and $g_t$ are PCF.
By Theorem \ref{DAOint}, for every  $t\in V(\C)$, $f_t$ and $g_t$ are intertwined.
By (3) and Theorem \ref{intertwined}, the images of $f_{t}$ and $g_{t}$ in $\sM_d$ are the same. This contradicts with (1), which concludes the proof.

\medskip

Recently DeMarco and Mavraki \cite{DeMarco2023} find a simplification of our proof of Theorem \ref{main}.   
Their  key observation is as follows: Our first key ingredient Theorem \ref{DAOint} proves more than what we need to show Theorem \ref{main}.
In Theorem \ref{DAOint}, we show that $f_t$ and $g_t$ are intertwined for {\bf every} $t\in V(\C).$ But we only need {\bf one} such parameter $t$ to prove Theorem \ref{main}. 
\subsection{Polynomial case.}
One can also consider the moduli space of polynomials and the multiplier spectrum morphism on it.  Let $\sP_d$ be the moduli space of degree $d$ polynomials, which  is the quotient space of the space of degree $d$ polynomials modulo the conjugation by affine maps. For every $d\geq 2$ and $n\geq 1$, we let $\tilde{\tau}_{d,n}$ (resp. $\tilde{\tau}_{d}$) be the restriction of $\tau_{d,n}$ (resp. $\tau_{d}$) on $\sP_d$. As the the moduli space of polynomials is a proper Zariski closed subset of $\sM_d$, Theorem \ref{main} does not implies the injectivity of $\tilde{\tau}_{d}$ directly. However, our proof of Theorem \ref{main} still works in the polynomial setting. 
\begin{thm}\label{thmpoly}For every $d\geq 2$,  $\tilde{\tau}_{d}$ is generically injective on $\sP_d$. 
\end{thm}
Huguin has an independent proof of Theorem \ref{thmpoly}, see \cite{huguin2024moduli}.
Indeed Huguin proved the generic injectivity of $\tilde{\tau}_{d,2}$. 
His method is completely different from ours in \cite{Ji2023b}, which is based on Fujimura's result \cite{Fujimura2007}  and on the computations  by Gorbovickis \cite{Gorbovickis2015}.

%

\subsection{Previous results.}
In \cite{Gorbovickis2015}, Gorbovickis showed that $\tau_{d,n}$ is generically quasi-finite for $d\geq 2$ and $n\geq 3.$
A recursive formula for the  upper bound of the topological degree of  $\tau_{d,n}$ was obtained by Schmitt in the preprint version \cite{Schmitt2016}.
An explicit upper bound of the topological degree of $\tau_{d,n}$ was obtained in Gotou's recent work \cite{Gotou2023}.

When $d=2$, Milnor \cite{Milnor1993} showed that  $\tau_{2,1}$ is in fact injective on $\sM_d(\C)$ (see also \cite[Theorem 2.45]{Silverman2012}). In particular Theorem \ref{main} holds when $d=2$.
When $d=3$,  Gotou showed that $\tau_{3,2}$ is generically injective (but not injective) \cite[Theorem 1.2]{Gotou2023}.
This was mentioned in \cite{Hutz2013a} which is the  errata for \cite{Hutz2013}.  
Previously there was no result about the generic injectivity of $\tau_d$ when $d\geq 4$.

For the polynomial case with $d\geq 2$, $\tilde{\tau}_{d,1}$ is generically quasi-finite (while $\tau_{d,1}$ is never generically quasi-finite, except when $d=2$).  Fujimura showed that $\deg (\tilde{\tau}_{d,1})=(d-2)!$ \cite{Fujimura2007}. The fiber structure of $\tilde{\tau}_{d,1}$ was studied by Sugiyama \cite{Sugiyama2017,Sugiyama2020,Sugiyama2023}.

\subsection{Further problems.}

\subsubsection{Injective locus of $\tau_d$.}
We have seen before that there are two mechanisms that can  produce rational maps with the same multiplier spectrum, one from Latt\`es maps, and another one from elementarily equivalent rational maps. Pakovich asked  \cite[Problem 3.1]{Pakovich2019} whether 
these are the only obstructions of the injectivity of $\tau_d$.  Motivated by Pakovich's question, Ji and me proposed the following conjecture.
\begin{con}\cite[Conjecture 1.7]{Ji2023b}\label{conequmulti}
	Let $f,g$ be rational maps of degree $d\geq 2$ such that the conjugacy classes of $f$ and $g$ are different. Assume that $\tau_d(f)=\tau_d(g)$, then one of the followings holds:
	\begin{enumerate}
		\item $f$ and $g$ are Latt\`es maps;
		\item $f$ is elementarily equivalent to $g$.
	\end{enumerate}
\end{con}


\subsubsection{Generic injectivity of multiplier spectrum of small periods.}
Our definition of $\tau_d$ requires the use of periodic points of periods not exceeding the number $m_d$, whose precise value is not effectively-known and is probably very large.  It is  interesting to know whether we can get generic injectivity using only periodic points of small periods.    By dimension counting, $\tau_{d,1}$ is never generically injective except when $d=2$.  Again by  dimension counting, it is very likely that  $\tau_{d,2}$ is generically injective. For this reason, Ji and me proposed the following conjecture:
\begin{con}\cite[Conjecture 1.8]{Ji2023b}
	For every $d\geq 2$,
	the morphism $$\tau_{d,2}:\sM_d(\C)\to \C^{N_1}/\Sigma_{N_1}\times \C^{N_2}/\Sigma_{N_2}$$ is generically injective. 
\end{con}

\subsubsection{Generic injectivity of the length spectrum.}
In \cite{Ji2023b}, Ji and me proposed the following conjecture, which is the analogy of Theorem \ref{main} for length spectrum. 
Note that $\sM_d$ is defined over $\Q$ (hence over $\R$).
\begin{con}\cite[Conjecture 1.9]{Ji2023b}\label{1.8}
	For every $d\geq 2$, there is a Zariski closed proper subset $E_d\subset \sM_d$ defined over $\R$, such that for every $[f]\notin E_d$, if there is $g\in\Rat_d(\C)$ such that $f$ and $g$ have the  same length spectrum, then  $g$ or $\overline{g}$ (the complex conjugation of $g$)  is conjugate to $f$.
\end{con}

By Theorem \ref{main}, Conjecture \ref{1.9} implies 
Conjecture \ref{1.8} directly.


\subsubsection{Multiplier spectrum for endomorphisms on $\P^N$.} 
For holomorphic endomorphisms on $\P^N$, $N\geq 2$, one can also construct the corresponding moduli space $\sM_d^N(\C)$ \cite{Silverman2012}. Moreover, the multiplier spectrum morphisms exist on  $\sM_d^N(\C)$ \cite{Silverman2012}. It is generically finite by Gauthier, Taflin and Vigny \cite[Corollary 2.4]{Gauthier2023}. It is of great interest to extend McMullen's multiplier spectrum rigidity (Theorem \ref{thmlength}) and the generic injectivity theorem (Theorem \ref{main}) to higher dimension.

We can ask same questions for other families of high dimensional dynamical systems such as H\'enon maps, automorphisms on K3 surfaces, ...

\section*{Acknowledgments.}
The author is supported by the NSFC Grant No.12271007.
The author would like to thank Zhuchao Ji, Xinyi Yuan and Charles Favre for their helpful comments, corrections, and suggestions.

\bibliographystyle{siamplain}
\bibliography{dd}
\end{document}